\documentclass[a4paper, DIV=calc, 11pt,  BCOR=0.5cm, twoside]{article} 

\usepackage[british]{babel} %Silbentrennung
\usepackage[utf8]{inputenc} %Inputencoding - für Sonderzeichen
\usepackage{microtype} %Randausgleich
\usepackage[
backend=biber,
sorting=nyt,
style=alphabetic,
maxbibnames=10
]{biblatex}
\renewbibmacro{in:}{}
\DeclareFieldFormat
  [article,inbook,incollection,inproceedings,patent,thesis,unpublished]
  {title}{#1\isdot}

\usepackage[toc]{appendix}

\usepackage[T1]{fontenc} %Fontencoding zum Schriften laden
\usepackage{lmodern} %Latin Modern
\usepackage{enumitem}
\setlist{  
	listparindent=\parindent,
	parsep=0pt,
}

\usepackage{dsfont, mathtools, amssymb, delarray, graphicx, nicefrac, import, cancel, mathrsfs, esint, framed, stmaryrd, pdfpages, amsmath}
\usepackage[pdfencoding=unicode, psdextra, colorlinks=false, pdfborder={0 0 0.7}]{hyperref}
\usepackage[style=english]{csquotes}

\usepackage[amsmath, thmmarks, hyperref]{ntheorem}

\theoremstyle{plain}
\theoremheaderfont{\normalfont\bfseries}
\theoremseparator{.}
\theorempreskip{\topsep}
\theorempostskip{\topsep}
\theoremnumbering{arabic}
\theoremsymbol{} %\NoEndMark
\newtheorem{Proposition}{Proposition}[section]
\newtheorem{Theorem}[Proposition]{Theorem}
\newtheorem*{Theorem-nn}{Theorem}
\newtheorem{Conjecture}[Proposition]{Conjecture}
\newtheorem{Question}[Proposition]{Question}
\newtheorem*{Conjecture-nn}{Conjecture}

\theoremstyle{plain}
\theoremheaderfont{\normalfont\bfseries}
\theoremseparator{.}
\theorempreskip{\topsep}
\theorempostskip{\topsep}
\theoremnumbering{arabic}
\newtheorem{Lemma}[Proposition]{Lemma}
\newtheorem*{Lemma-nn}{Lemma}
\newtheorem{Definition}[Proposition]{Definition}
\newtheorem*{Definition-nn}{Definition}
\newtheorem{Corollary}[Proposition]{Corollary}

\theoremstyle{plain}
\theoremheaderfont{\normalfont\bfseries}
\theorembodyfont{\upshape}
\theoremseparator{.}
\theorempreskip{\topsep}
\theorempostskip{\topsep}
\theoremnumbering{arabic}
\theoremsymbol{}

\newtheorem*{Reminder-nn}{Reminder}
\newtheorem{Remark}[Proposition]{Remark}
\newtheorem*{Remark-nn}{Remark}

\newtheorem*{Example-nn}{Example}

\newtheorem*{Motivation-nn}{Motivation}
\newtheorem{Construction}[Proposition]{Construction}

\newtheorem*{Notation-nn}{Notation}

\theoremstyle{plain}
\theoremheaderfont{\normalfont\bfseries}
\theoremseparator{}
\theorempreskip{\topsep}
\theorempostskip{\topsep}
\theoremnumbering{arabic}
\theoremsymbol{}

\theoremstyle{plain}							
\theoremheaderfont{\itshape}
\theorembodyfont{\upshape}
\theoremseparator{.}
\theorempreskip{\topsep}
\theorempostskip{\topsep}
\theoremnumbering{arabic}
\theoremsymbol{\ensuremath{\square}}
\newtheorem*{Proof}{Proof}

\theoremstyle{plain}							
\theoremheaderfont{\itshape}
\theorembodyfont{\upshape}
\theoremseparator{.}
\theorempreskip{\topsep}
\theorempostskip{\topsep}
\theoremnumbering{arabic}
\theoremsymbol{\ensuremath{\square}}

\theoremstyle{plain}							
\theoremheaderfont{\normalfont\bfseries}
\theoremseparator{}
\theorempreskip{\topsep}
\theorempostskip{\topsep}
\theoremnumbering{arabic}
\theoremsymbol{}
\newtheorem*{Blank}{}

\usepackage[capitalize,nameinlink]{cleveref}
\crefname{Lemma}{Lemma}{Lemmata}
\crefname{Definition}{Definition}{Definitions}
\crefname{Example}{Example}{Examples}
\crefname{Theorem}{Theorem}{Theorems}
\crefname{Corollary}{Corollary}{Corollaries}
\crefname{Proposition}{Proposition}{Propositions}
\crefname{Equation}{Eq.}{Eq.}
\crefname{IEEEeqnarray}{Eq.}{Eq.}
\crefname{Remark}{Remark}{Remarks}
\crefname{Reminder}{Reminder}{Reminders}
\crefname{figure}{Figure}{Figures}
\crefname{chapter}{Chapter}{Chapters}
\crefname{Construction}{Construction}{Constructions}
\crefname{align}{Eq.}{Eq.}
\crefname{align*}{Eq.}{Eq.}
\crefname{Conjecture}{Conjecture}{Conjectures}
\crefname{Question}{Question}{Questions}

\let\existstemp\exists
\let\foralltemp\forall
\let\nexiststemp\nexists

\renewcommand*{\exists}{\mkern3mu\existstemp\mkern3mu}
\renewcommand*{\forall}{\mkern3mu\foralltemp\mkern3mu}
\renewcommand*{\nexists}{\mkern3mu\nexiststemp\mkern3mu}

\usepackage{IEEEtrantools}

\makeatletter
\newcommand{\oset}[3][0ex]{%
	\mathrel{\mathop{#3}\limits^{
			\vbox to#1{\kern-2\ex@
				\hbox{$\scriptstyle#2$}\vss}}}}
\makeatother

\renewcommand{\imath}{\mathrm{i}\mkern1mu}

\usepackage{color}

\makeatletter
\newcommand{\dashedrightarrow}[1][2pt]{%
	\settowidth{\@tempdima}{$\rightarrow$}\rightarrow% typeset arrow
	\makebox[-\@tempdima]{\hskip-2ex\color{white}\rule[0.5ex]{#1}{1pt}}% typeset overlay
	\phantom{\rightarrow}% advance appropriate horizontal distance
}
\makeatother

\DeclareMathOperator{\Hom}{Hom}
\DeclareMathOperator{\Aut}{Aut}
\DeclareMathOperator{\Ext1}{Ext^1_\mathcal{O}}

\DeclareMathOperator{\id}{id}

\DeclareMathOperator{\Alb}{Alb}
\DeclareMathOperator{\vol}{vol}

\usepackage{fancyhdr}
\title{Canonical extensions of manifolds with nef tangent bundle}
\author{Niklas Müller}
\date{}

\begin{document}
	\maketitle
	
	\begin{abstract}
	
	\noindent
		To any compact Kähler manifold $(X, \omega)$ one may associate a bundle of affine spaces $Z_X\rightarrow X$ called a \emph{canonical extension} of $X$. In this paper we prove that if the tangent bundle of $X$ is nef, then the total space $Z_X$ is a Stein manifold. This partially answers a question raised by Greb-Wong of whether these two properties are actually equivalent. We also complement some known results for surfaces in the converse direction.
	\end{abstract}

	\setcounter{section}{-1}
	\section{Introduction}
	Given a Kähler manifold $(X, \omega)$ one may define in a natural way a bundle $p\colon Z_X\rightarrow X$ of affine spaces called a \emph{canonical extension} of $X$. One possible way to define $Z_X$ is as the universal complex manifold on which the cohomology class $[p^*\omega]=0$ vanishes. 
	
	Canonical extensions were introduced by \cite{donaldson_SmoothnessOfMongeAmpereFlow} to prove regularity properties of solutions to the Monge-Ampère equation. They have subsequently also seen some uses related to K-stability and the existence of Kähler-Einstein metrics on Fano manifolds, see for example \cite{tian_KStabilityStrongly} or \cite{gkp_CanExtandKSTability}. Recently, in \cite{greb_CanonicalExtensions}, the following question was posed which suggests another point of view on canonical extensions:
	\begin{Question}
		Let $X$ be a compact Kähler manifold. Is it true, that the tangent bundle of $X$ is nef if and only if some (resp.\ any) canonical extension of $X$ is a Stein manifold?\label{IntroQ:GrebWong}
	\end{Question}
	The structure of compact Kähler manifolds possessing a nef tangent bundle is well-understood and by now classical. However, specifically in the Fano case some very interesting questions such as the conjecture of Campana-Peternell remain open. Thus, \cref{IntroQ:GrebWong} is interesting as it suggests a possibly more geometric point of view on these problems. In \textbf{\cref{sec:canExtGeneralManifoldsWithNefTangent}} we give the following partial answer to \cref{IntroQ:GrebWong}:
	\begin{Theorem}
		Let $X$ be a compact K\"ahler manifold with nef tangent bundle. If the (weak) Campana-Peternell conjecture \cref{con:Campana-Peternell} holds true\footnote{Note added in proof: Recently, Conjecture~\ref{con:Campana-Peternell} has been proved by Wang \cite{wang_WeakCampanaPeternellConjecture}. In particular, the conclusion in Theorem~\ref{IntroCanExtStein} holds unconditionally.}
        then any canonical extension $Z_X$ of $X$ is a Stein manifold.\label{IntroCanExtStein}
	\end{Theorem}
    \cref{IntroCanExtStein} was previously only known to hold in the special cases of complex tori by \cite[Proposition 2.13]{greb_CanonicalExtensions} and for Fano manifolds with big tangent bundle, see \cite[Theorem 1.2]{HP_Stein_complements}.
	
	In the converse direction of what can be said about manifolds admitting a canonical extension which is Stein, little is known. In fact, even for projective surfaces \cref{IntroQ:GrebWong} is not completely settled yet, although it is known to hold in most cases by the work of \cite[Theorem 1.13]{HP_Stein_complements}. In \textbf{\cref{sec:CanExtSurfaces}} we partially complement their results by treating also the case of ruled surfaces over curves of higher genus:
	\begin{Lemma}
		Let $X = \mathds{P}(\mathcal{E}) \rightarrow C$ be a ruled surface over a curve of genus $g(C)\geq 2$ defined by a semi-stable vector bundle $\mathcal{E}$. Then no canonical extension of $X$ is Stein.\label{Intro:RuledSurfaces}
	\end{Lemma}
	This only leaves to consider the case of unstably ruled surfaces over elliptic curves.

	\subsection*{Acknowledgements}
	I would like to express my gratitude towards Daniel Greb, my supervisor, who suggested the main question that is answered in this publication as a master's thesis topic. Moreover, I am incredibly thankful for his very valuable advice.

	\section{Canonical extensions of complex manifolds}
	\noindent
	In this section we describe a general approach to constructing bundles of affine spaces over complex manifolds and define the canonical extension. Let $X$ be a complex manifold and fix a holomorphic vector bundle $\mathcal{E}$ on $X$ and a cohomology class $a\in H^1(X, \mathcal{E}) = \textmd{Ext}^1_{\mathcal{O}}(\mathcal{O}_X, \mathcal{E})$. Denote by $0 \rightarrow \mathcal{E} \rightarrow \mathcal{V}_a \overset{p}{\rightarrow} \mathcal{O}_X \rightarrow 0$ the extension corresponding to $a$. Below we describe three equivalent ways of constructing affine bundles over $X$ from the data $(\mathcal{E}, a)$:
	\begin{Construction}
		(\textbf{as torsors})
		
		\noindent
		Consider the sub sheaf $\mathcal{Z}_a := p^{-1}(1) \subsetneq \mathcal{V}_a$ of sections mapping to the constant function $1$. Note that $\mathcal{Z}_a$ comes with a natural action of $\mathcal{E}$ by translations making $\mathcal{Z}_a$ into an \emph{affine bundle} in the following sense: The underlying total space $Z_a := |\mathcal{Z}_a|\rightarrow X$ is a fibre bundle over $X$ and the fibre $Z_a|_x$ over any point $x$ is in a natural way an affine vector space with group of translations $\mathcal{E}|_x$. We call $Z_a \rightarrow X$ an \emph{extension} of $X$ \emph{modelled on the vector bundle $\mathcal{E}$}. We may also denote $Z_{\mathcal{E}, a}$ if we want to emphasise the role of $\mathcal{E}$. 

        Equivalently, $Z_a = |p|^{-1}(X\times \{ 1 \})$, where $|p|\colon |\mathcal{V}_a|\rightarrow |\mathcal{O}_X| = X\times\mathds{C}$ denotes the holomorphic map between the underlying total spaces of the bundles $\mathcal{V}_a$, $\mathcal{O}_X$. This is the definition of $Z_a$ used in \cite{greb_CanonicalExtensions}.\label{con:CanonicalExtensions}
	\end{Construction}

	\begin{Construction}
		(\textbf{as complements of a hypersurface})	
		
		\noindent	
		A second, possibly more geometric construction of $Z_a$ is as follows: Consider dually the short exact sequence
		\begin{align*}
		0 \rightarrow \mathcal{O}_X \rightarrow (\mathcal{V}_a)^* \rightarrow \mathcal{E}^*  \rightarrow 0
		\end{align*}
		which defines an embedding $\mathds{P}(\mathcal{E}^*) \hookrightarrow \mathds{P}(\mathcal{V}_a^*)$. Here, throughout this paper we will always use the convention that $\mathds{P}(\mathcal{E})$ denotes the projective bundle of linear \emph{hyperplanes} in $\mathcal{E}$. 

        It is then clear from the construction that here exists a natural identification 
        \begin{align*}
            Z_a = \mathds{P}(\mathcal{V}_a^*)\setminus \mathds{P}(\mathcal{E}^*).
        \end{align*}
		Note that by construction $\mathds{P}(\mathcal{E}^*)$ is embedded as a smooth hypersurface in the linear series of $\mathcal{O}_{\mathds{P}(\mathcal{V}_a^*)}(1)$. In particular, its normal bundle is given by $\mathcal{N}_{\mathds{P}(\mathcal{E}^*) / \mathds{P}(\mathcal{V}_a^*)}
			= \mathcal{O}_{\mathds{P}(\mathcal{E}^*)}(1)$. 
			This is the preferred point of view in \cite{HP_Stein_complements}.
			
	\end{Construction}
	
	\begin{Construction}
		(\textbf{via a universal property})
		
		\noindent
		Finally, $Z_a \overset{p}{\rightarrow} X$ enjoys the following universal property, which of course determines it uniquely: Let $h\colon Y \rightarrow X$ be any holomorphic map from a complex manifold such that $h^*a = 0\in H^1(Y, f^*\mathcal{E})$. Then $h$ factors uniquely, up to translation by an element of $H^0(Y, f^*\mathcal{E})$, through $Z_a\overset{p}{\rightarrow} X$. In this sense, $Z_a\rightarrow X$ is the universal manifold on which the cohomology class $a$ vanishes. A more precise version of this statement may be found in \cite[Lemma 1.16.(c)]{greb_CanonicalExtensions}.\label{universal_property_can_extensions}
	\end{Construction}

	\begin{Definition}
		Let $(X, \omega)$ be a complex Kähler manifold. Then $\omega$ is a $\bar{\partial}$-closed form and hence defines a cohomology class $[\omega] \in H^1(X, \Omega_X^1)$. The associated extension $Z_X := Z_{[\omega]}$ is called (a) \emph{canonical extension} of $X$.
	\end{Definition}
	Above, we have seen three equivalent constructions for $Z_{[\omega]}$:
	\begin{itemize}
		\item[(1)] As a bundle of affine spaces over $X$ modelled on the cotangent bundle $\Omega^1_X$.
		\item[(2)] As the complement $Z_{[\omega]} = \mathds{P}(\mathcal{V}^*)\setminus \mathds{P}(\mathcal{T}_X)$ of the smooth hypersurface $\mathds{P}(\mathcal{T}_X)$ whose normal bundle is given by $\mathcal{N}_{\mathds{P}(\mathcal{T}_X)/ \mathds{P}(\mathcal{V}^*)} = \mathcal{O}_{\mathds{P}(\mathcal{T}_X)}(1)$.
		\item[(3)] As the universal manifold on which the cohomology class $[\omega]$ vanishes.
	\end{itemize}
	The following conjecture arose out of the work of \cite{greb_CanonicalExtensions} and \cite{HP_Stein_complements} on canonical extensions:
	\begin{Conjecture}\emph{(Greb-Wong, Höring-Peternell)}
		
		\noindent
		Let $X$ be a compact Kähler manifold. Then the tangent bundle $\mathcal{T}_X$ is nef if and only if some canonical extension $Z_X$ of $X$ is Stein.\label{con:GrebWongCanExt}
	\end{Conjecture}
    \cref{con:GrebWongCanExt} was confirmed for Kähler manifolds of non-negative holomorphic bisectional curvature, e.g.\ complex tori and flag manifolds by \cite{greb_CanonicalExtensions}. Moreover, assuming that $X$ is a Fano manifold with big and nef tangent bundle it was shown in \cite[Theorem 1.2]{HP_Stein_complements} that any canonical extension of $X$ must be affine and, hence, Stein. Below, we will combine both cases to prove \cref{IntroCanExtStein}.
    \begin{Remark}
        Similarly, Greb-Wong and Höring-Peternell conjectured the tangent bundle of a smooth projective variety $X$ should be big and nef if and only if some canonical extension of $X$ is affine. Combining results of \cite[Corollary 4.4]{greb_CanonicalExtensions} and \cite[Theorem 1.2]{HP_Stein_complements} this second conjecture may be reduced to the first one, \cref{con:GrebWongCanExt}.
    \end{Remark}

	\section{Structure theory of manifolds with nef tangent bundle}
	\label{sec:StructureManifoldsWithNefTangent}

    The following result will be a key ingredient in the proof of \cref{IntroCanExtStein}; it summarises the successive work of \cite{campanaPeternell_Conjecture}, \cite{DPS_ManifoldsWithNefRicciClass}, \cite{DPS_ManifoldsWithNefTangentBundle} and \cite{cao_PhDThesis}:
	
	\begin{Theorem}\emph{(Cao, Demailly-Peternell-Schneider)}
		
		\noindent
		Let $X$ be a compact Kähler manifold possessing a nef tangent bundle. There exists a finite étale cover $X' \rightarrow X$ such that the Albanese map $\alpha\colon X' \rightarrow \Alb(X')$ is a locally constant holomorphic fibre bundle. The typical fibre is a Fano manifold with nef tangent bundle. \label{flatness_of_Albanese}
	\end{Theorem}
	Here, a fibre bundle $\alpha\colon X \rightarrow T$ is said to be \emph{locally constant} if its transition functions may be chosen to be locally constant. This is equivalent to the existence of a group homomorphism $\rho\colon \pi_1(T) \rightarrow \Aut(F)$ such that $(\widetilde{T}\times F)/\pi_1(T) \cong X$ as fibre bundles over $T$. Note that in the latter case $pr_{\widetilde{T}}^*\mathcal{T}_{\widetilde{T}}$ descends to a holomorphic vector bundle on $X$ providing a global holomorphic splitting for the short exact sequence
    \begin{align}
        0 \rightarrow \mathcal{T}_{X/T} \rightarrow \mathcal{T}_{X} \rightarrow \alpha^*\mathcal{T}_{T} \rightarrow 0. \label{eqs:LocallyConstantFibreBundle}
    \end{align} 
		
	The following famous conjecture claims that in the situation of \cref{flatness_of_Albanese} much more can be said about the fibre of $\alpha$:
	\begin{Conjecture}\emph{(Campana-Peternell, \cite{campanaPeternell_Conjecture}, weakend form)}
    \label{con:Campana-Peternell}
		
		\noindent 
        If the tangent bundle of a Fano manifold is nef then it must also be big\footnote{Note added in proof: As stated already in the introduction, Conjecture~\ref{con:Campana-Peternell} has been recently proved by Wang \cite{wang_WeakCampanaPeternellConjecture}.}.
	\end{Conjecture}
    In fact, the original formulation of \cref{con:Campana-Peternell} is stronger and predicts that any Fano manifold with nef tangent bundle should even be homogeneous. That the tangent bundle of a homogeneous Fano manifold must be big is a classical fact; see e.g.\ \cite[Corollary 4.4]{greb_CanonicalExtensions} for a proof using canonical extensions or, alternatively, \cite[Corollary 1.3]{hsiao_BignessTangentBundleFlagVariety}. 
	
	The conjecture of Campana and Peternell has seen attention by quite a number of authors and is by now verified for manifolds of dimension at most five by \cite{kanemitsu_Fano5foldsNefTangentBundle}, see also the introduction thereof for a short summary of contributions to this problem. In full generality however even its weaker form \cref{con:Campana-Peternell} is still completely open.

	\section{Canonical extensions of manifolds with nef tangent bundle}
	\label{sec:canExtGeneralManifoldsWithNefTangent}
	
	In this section we prove \cref{IntroCanExtStein}. Let us fix a compact Kähler manifold $(X, \omega_X)$.

	\begin{Proposition}
		Assume that the Albanese morphism $\alpha\colon X \rightarrow \Alb(X) =:T$ is a locally constant holomorphic fibre bundle. Then there exists a natural isomorphism of affine bundles 
		\begin{align*}
		Z_{X, [\omega_X]} 
		\cong Z_{\Omega^1_{X/T}, [\omega_{X/T}]} \times_X Z_{\alpha^*\Omega^1_T, a_T}.
		\end{align*}
		Here, $[\omega_{X/T}]$ denotes the image of $[\omega_X]$ under the natural map $H^1(X, \Omega_X^1) \rightarrow H^1(X, \Omega_{X/T}^1)$.\label{splitting_canonical_extension_2}
	\end{Proposition}
	In the above statement we leave the extension class $a_T$ ambiguous on purpose. It will be described more explicitly in \cref{determination_of_pushdown_class} below.
	\begin{Proof}[of \Cref{splitting_canonical_extension_2}]
		As $\alpha$ is a locally constant fibraton we have $\mathcal{T}_{X} \cong \mathcal{T}_{X/T} \oplus \alpha^*\mathcal{T}_{T}$ (see  the discussion around \cref{eqs:LocallyConstantFibreBundle}). Hence, it follows from \cite[Lemma 5.5]{HP_Stein_complements} that
		\begin{align*}
		Z_{X, [\omega_X]} := Z_{\Omega^1_X, [\omega_X]}
		\cong Z_{\Omega^1_{X/T}, [\omega_{X/T}]} \times_X Z_{\alpha^*\Omega^1_T, a_T},
		\end{align*}
        where $[\omega_X] = [\omega_{X/T}] \oplus a_T \in \Ext1(\mathcal{O}_X, \Omega_X^1) \cong \Ext1(\mathcal{O}_X, \Omega_{X/T}^1) \oplus \Ext1(\mathcal{O}_X, \alpha^*\Omega^1_T )$ is the induced decomposition. In other words, $[\omega_{X/T}]$ is the image of $[\omega_X]$ under the natural map
		\begin{align*}
		\Ext1\left(\mathcal{O}_X, \Omega_X^1\right) \rightarrow \Ext1\left(\mathcal{O}_X, \Omega_{X/T}^1\right).
		\end{align*}
		Modulo the identification $\Ext1(\mathcal{O}_X, - ) = H^1(X, -)$ this is the proclaimed class.
	\end{Proof}
	Our next goal is to give an explicit description of the cohomology class $a_T$ in \cref{splitting_canonical_extension_2}. To this end, let $f\colon X \rightarrow T$ be a submersion of relative dimension $m$. Let us denote by $F_t$ the fibres of $f$. Then the function
		\begin{align*}
		\vol(F_t, \omega_X|_{F_t}) := \frac{1}{m!} \int_{F_t} \left(\omega_X|_{F_t}\right)^m = \frac{1}{m!} \hspace{0.1cm} f_* \left(\omega_X^m\right)\big|_t
		\end{align*}
		is constant. Here, $f_*$ denotes the \emph{integration along the fibres} and the constancy of $f_* (\omega_X^m)$ is clear as $f_*$ commutes with the exterior derivative and as $\omega_X$ is $d$-closed.
	
	\begin{Proposition}
    If any fibre $F_t$ of the submersion $f\colon X \rightarrow T$ is Fano then the composition
		\begin{align*}
		P\colon H^q\left(X, f^*\Omega_{T}^p\right) \overset{i_*}{\longrightarrow} H^q\left(X, \Omega_X^p\right) \xrightarrow{\wedge\frac{\omega^m}{m!}} H^{q+m}\left(X, \Omega_X^{p+m}\right) \overset{f_*}{\longrightarrow} H^q\left(T, \Omega_T^p\right)
		\end{align*}
		is an isomorphism for all $p,q$. In fact, the inverse is given up to a scalar factor by the natural map $f^*\colon H^q\left(T, \Omega_T^p\right) \rightarrow H^q\left(X, f^*\Omega_{T}^p\right)$.\label{integration_along_fibres}
	\end{Proposition}
	
	\begin{Proof}
    First, $R^jf_*f^*\Omega_T^p = \Omega_{T}^p\otimes R^jf_*\mathcal{O}_X = \Omega_{T}^p\otimes R^jf_*\mathcal{O}_X(-K_X+K_X) = 0$ for all $j>0$ due to the Kodaira vanishing theorem. Thus, it follows from the Leray spectral sequence that $f^*$ is an  isomorphism. Below, we will prove using Dolbeaut representatives that up to a scalar factor its inverse is given by $P$: 
    Fix any integers $p,q$ and any closed differentiable $(p,q)$-form $\eta$ on $T$. We compute
		\begin{align}
		P(f^*([\eta])) 
		&=: \frac{1}{m!} \hspace{0.1cm} \left[f_*\big(f^*\eta\wedge \omega_X^m \big)\right] \nonumber =
		\frac{1}{m!} \hspace{0.1cm} \left[ \eta \wedge f_*(\omega_X^m)\right] \nonumber =: [\eta] \cdot \vol(F) \label{eqs:ClassDefiningCanExt}
		\end{align}
		so that $P\circ f^* = \vol(F)\cdot \id$. This concludes the proof.
	\end{Proof}
	\begin{Proposition}
		Assume that any fibre of $f\colon X \rightarrow T$ is a Fano manifold and that the natural short exact sequence $0 \rightarrow f^*\Omega_{T}^1 \rightarrow \Omega_X^1 \rightarrow \Omega_{X/T}^1 \rightarrow 0$ admits a splitting $s\colon \Omega_X^1\rightarrow f^*\Omega_{T}^1$. Let
		\begin{align*}
		[\omega_X] = [\omega_{X/T}] + a_T \in H^1\left(X, \Omega_X^1\right) = H^1\left(X, \Omega_{X/T}^1\right) \oplus H^1\left(X, f^*\Omega_{T}^1\right)
		\end{align*}
		be the induced decomposition so that $a_T = H^1(s)([\omega_X])$ and let $\omega_T := f_*(\omega_X^{m+1})$ denote the Kähler form on $T$ obtained from $\omega_X$ by integration along the fibres. Then
		\begin{align}
		a_T = \frac{1}{(m+1)!\cdot \vol(F)} \cdot [f^*\omega_T] \in H^1\left(X, \Omega_X^1\right).\label{eqs:ClassDefiningCanExt_2}
		\end{align}\label{determination_of_pushdown_class}
	\end{Proposition}
	\begin{Corollary}
		Let $(X, \omega_X)$ be a compact Kähler manifold with nef tangent bundle. If the Albanese $\alpha\colon X \rightarrow \Alb(X) =:T$ is a locally constant holomorphic fibre bundle with Fano manifolds as fibres then there exists a natural isomorphism of affine bundles 
		\begin{align*}
		Z_{X, [\omega_X]} 
		\cong Z_{\Omega^1_{X/T}, [\omega_{X/T}]} \times_X Z_{\alpha^*\Omega^1_T, [\alpha^*\omega_T]}
		\cong Z_{\Omega^1_{X/T}, [\omega_{X/T}]} \times_T Z_{T, [\omega_T]}.
		\end{align*}\label{fibreBundleStructureOnZX}
	\end{Corollary}
	\begin{Proof}
		Since $\alpha$ is locally constant the sequence $0 \rightarrow \alpha^*\Omega_{T}^1 \rightarrow \Omega_X^1 \rightarrow \Omega_{X/T}^1 \rightarrow 0$ splits. Moreover, according to \cref{determination_of_pushdown_class} the induced decomposition of $[\omega_X]$ is given by
		\begin{align*}
		\left[\omega_X\right] = \left[\omega_{X/T}\right] + \lambda \cdot \left[\alpha^*\omega_T\right] \in  \Ext1\left(\mathcal{O}_X, \Omega_X^1\right) = \Ext1\left(\mathcal{O}_X, \Omega_{X/T}^1\right)\oplus \Ext1\left(\mathcal{O}_X, \alpha^*\Omega_{T}^1\right),
		\end{align*} 
		where $\lambda := \frac{1}{(m+1)!\cdot \vol(F)}$ is some constant. In effect, the proof of \cref{splitting_canonical_extension_2} shows that
		\begin{align*}
		Z_{X, [\omega_X]} 
		\cong Z_{\Omega^1_{X/T}, [\omega_{X/T}]} \times_X Z_{ \alpha^*\Omega^1_T, \lambda \cdot [\alpha^*\omega_T]}.
		\end{align*}
		Since extensions only depend on their defining cohomology class up to scaling by \cite[Remark 2.4]{greb_CanonicalExtensions} it follows that
		\begin{align*}
		Z_{X, [\omega_X]} 
		\cong Z_{\Omega^1_{X/T}, [\omega_{X/T}]} \times_X Z_{ \alpha^*\Omega^1_T, [\alpha^*\omega_T]}
		\cong Z_{\Omega^1_{X/T}, [\omega_{X/T}]} \times_T Z_{\Omega^1_T, [\omega_T]}.
		\end{align*}
		Here in the last step we used the functoriality of extensions, see \cite[Lemma 1.16(b)]{greb_CanonicalExtensions}. 
	\end{Proof}

	\begin{Proof}[of \cref{determination_of_pushdown_class}]
	
		\noindent
		We will verify \cref{eqs:ClassDefiningCanExt_2} by an explicit calculation using Dolbeaut representatives. To this end, recall that $s\colon \Omega_X^1\rightarrow f^*\Omega_{T}^1$ induces maps of sections $s^{(0,1)}\colon\mathcal{A}^{0,1}( \Omega_X^1) \rightarrow \mathcal{A}^{0,1}(f^*\Omega_{T}^1)$ and the class 
		\begin{align}
		i_*(a_T) = i_*\left(H^1(s)\big([\omega_X]\big)\right) \in H^1\left(X, f^*\Omega_T^1\right) \overset{i_*}{\hookrightarrow} H^1\left(X, \Omega^1_X\right) \label{eqs:ClassDefiningCanExt_2.5}
		\end{align}
		is represented by the form $\zeta := i_*(s^{(0,1)}(\omega_X))$. Below, we will show that
		\begin{align}
		f_*(\zeta\wedge\omega_X^m) =  \frac{f_*(\omega_X^{m+1})}{m+1}.  \label{eqs:ClassDefiningCanExt_3}
		\end{align}
		This will immediately yield the result because, assuming \cref{eqs:ClassDefiningCanExt_3} and using \cref{integration_along_fibres}, we compute
		\begin{align}
		i_*(a_T) &=: \left[\zeta\right] 
		= \frac{1}{\vol(F)} \cdot i_*\left[f^*f_*\left(\zeta\wedge\frac{\omega_X^m}{m!}\right)\right] \overset{\textmd{\cref{eqs:ClassDefiningCanExt_3}}}{=\joinrel=\joinrel=\joinrel=} \frac{1}{\vol(F)} \cdot \frac{1}{(m+1)!} \cdot i_*\left[f^*f_*\big(\omega_X^{m+1}\big)\right]\nonumber \\
		&\qquad =: \frac{1}{\vol(F) \cdot (m+1)!} \cdot i_*\left[f^*\omega_T\right].\label{eqs:ClassDefiningCanExt_3.5}
		\end{align}
		As $i_*$ is injective by \cref{integration_along_fibres} this is the equation to prove. In conclusion, it remains to verify \cref{eqs:ClassDefiningCanExt_3}. To this end, fix a point $t\in T $ and vectors $v \in T_t^{(1,0)}T$, $w \in T_t^{(0,1)}T$. Let $\widetilde{V} := s^*(v), \widetilde{W} := s^*(w)$ be the differentiable vector fields along $F_t$ induced by the dual splitting $s^*\colon f^*\mathcal{T}_T\hookrightarrow \mathcal{T}_X$. Then $\widetilde{V}, \widetilde{W}$ are of type $(1,0)$ (respectively $(0,1)$) and lift $v, w$:
		\begin{align*}
		df(\widetilde{V}|_x) = v, \quad df(\widetilde{W}|_x) = w, \quad \forall x\in F_t.
		\end{align*}
		By definition, we have the identities
		\begin{align}
		\left(f_*\left( \zeta\wedge\omega_X^m \right)\right) (v, w) &= \int_{F_t} \iota_{\widetilde{V}, \widetilde{W}}\left(\zeta\wedge\omega_X^m\right), \label{eqs:ClassDefiningCanExt_4}\\
		\left(f_*\omega_X^{m+1}\right) (v, w) &= \int_{F_t} \iota_{\widetilde{V}, \widetilde{W}}\left( \omega_X^{m+1} \right)\label{eqs:ClassDefiningCanExt_5}
		\end{align}
		and we need to prove the equality of both expressions (modulo a scalar factor). Clearly it suffices to prove point-wise equality of the integrands as differential forms and this is what we will do: Fix a point $x\in F_t$ and denote $\widetilde{v} := \widetilde{V}|_x, \widetilde{w} := \widetilde{W}|_x$. 
		
		\begin{Blank}
			\emph{Step 1: For all tangent vectors $v'\in T^{1,0}_xX, w'\in T^{0,1}_xX$ it holds that}
			\begin{align*}
			\zeta(v', w') \overset{\textmd{\cref{eqs:ClassDefiningCanExt_2.5}}}{:=\joinrel=\joinrel=\joinrel=\joinrel=} 
			i_*\left(s^{(0,1)}\left(\omega_X\right)\right)(v', w') 
			= \omega_X\left(s^*\left(df(v')\right), w'\right).
			\end{align*}
		\end{Blank}
		Indeed, if more generally $\phi\colon \mathcal{E} \rightarrow \mathcal{F}$ is any morphism between holomorphic vector bundles, then the induced map $\phi^{(0,1)}\colon\mathcal{A}^{0,1}(\mathcal{E}) \rightarrow \mathcal{A}^{0,1}(\mathcal{F})$	is determined by the rule $\phi^{(0,1)}(\sigma\otimes d\bar{z}) = \phi(\sigma) \otimes d\bar{z}$. Accordingly, if $(z^j)$ are some local coordinates centred at $x\in F_t$ and if with respect to these coordinates $\omega_X =  \sum h_{k, \ell}\hspace{0.1cm} dz^k\wedge d\bar{z}^{\ell}$, then $s^{(0,1)}(\omega_X)$ is locally given by the expression
		\begin{align*}
		s^{(0,1)}(\omega_X) = s^{(0,1)}\left(\sum h_{k, \ell}\hspace{0.1cm} dz^k\wedge d\bar{z}^{\ell}\right) = \sum h_{k, \ell} \hspace{0.1cm} s\left(dz^k\right) \otimes d\bar{z}^{\ell}.
		\end{align*}
		Similarly, $i_*\colon \mathcal{A}^{0,1}(f^*\Omega_T^1) \hookrightarrow \mathcal{A}^{0,1}(\Omega_X^1)$ is by construction the map induced by the bundle morphism $(df)^*\colon f^*\Omega^1_T \hookrightarrow \Omega_X^1$. In other words,
		\begin{align*}
		i_*\left(s^{(0,1)}\left(\omega_X\right)\right)(v', w')
		: &= \left( \sum h_{k, \ell} \hspace{0.1cm} df^*\big(s\left(dz^k\right)\big) \otimes d\bar{z}^{\ell} \right)(v', w')\\
		&= \sum h_{k, \ell } \hspace{0.1cm} \big((df^*\circ s)(dz^k)\big)(v') \otimes d\bar{z}^{\ell}(w') \\
		&= \sum h_{k, \ell } \hspace{0.1cm} dz^k\big(s^*(df(v'))\big) \otimes d\bar{z}^{\ell}(w') \\
		&= \left( \sum h_{k, \ell} \hspace{0.1cm} dz^k \otimes d\bar{z}^{\ell} \right)\left(s^*(df(v')), w'\right)
		= \omega_X\left(s^*(df(v')), w'\right).
		\end{align*}

		\begin{Blank}
			\emph{Step 2: We have} $\iota_{\widetilde{v}, \widetilde{w}}(\zeta\wedge \omega_X^m)\big|_{F_t}  = \big(\omega_X(\widetilde{v}, \widetilde{w}) \cdot \omega_X^m - \iota_{\widetilde{v}}(\omega_X) \wedge \iota_{\widetilde{w}}(\omega_X^m) \big)\big|_{F_t}.$
		\end{Blank}
		Using some elementary formulae from multi-linear algebra we compute
		\begin{align}
		\iota_{\widetilde{w}}\iota_{\widetilde{v}}(\zeta\wedge \omega_X^m)
		&= \iota_{\widetilde{w}}\left( \iota_{\widetilde{v}}(\zeta)\wedge \omega_X^m + (-1)^2 \hspace{0.1cm} \zeta\wedge \iota_{\widetilde{v}}(\omega_X^m) \right)\nonumber\\
		&= \zeta(\widetilde{v}, \widetilde{w}) \cdot \omega_X^m + (-1) \hspace{0.1cm} \iota_{\widetilde{v}}(\zeta)\wedge \iota_{\widetilde{w}}(\omega_X^m)\nonumber\\
		&\qquad + (-1)^2\hspace{0.1cm} \iota_{\widetilde{w}}(\zeta)\wedge \iota_{\widetilde{v}}(\omega_X^m)
		+ (-1)^4\hspace{0.1cm} \zeta \wedge \iota_{\widetilde{v}, \widetilde{w}}(\omega_X^m).\label{eqs:ClassDefiningCanExt_10}
		\end{align}
		Now, according to \emph{Step 1} it holds that 
        \begin{align}
            \zeta(v', -) = \omega_X(s^*(df(v')), -), \quad \forall v'\in T^{0,1}_xX.\label{eqs:ClassDefiningCanExt_11} 
        \end{align}
        In particular, if $v'$ is tangent along the fibres, then $df(v') = 0$ and so $\iota_{v'}\zeta = 0$. Thus,
		\begin{align}
		\iota_{\widetilde{w}}(\zeta)|_{F_t} = \zeta|_{F_t} = 0.\label{eqs:ClassDefiningCanExt_12}
		\end{align}
		On the other hand, consider the case $v' = \widetilde{v}$ in \cref{eqs:ClassDefiningCanExt_11} above. Then $
		s^*(df(\widetilde{v})) = s^*(v) = \widetilde{v}$
		by definition of $\widetilde{v}$. In view of \cref{eqs:ClassDefiningCanExt_11} this implies that
		\begin{align}
		   \iota_{\widetilde{v}}(\zeta) = \iota_{\widetilde{v}}(\omega_X).\label{eqs:ClassDefiningCanExt_13}
		\end{align}
		Substituting the terms in \cref{eqs:ClassDefiningCanExt_10} above using \cref{eqs:ClassDefiningCanExt_12} and \cref{eqs:ClassDefiningCanExt_13} we find
		\begin{align*}
		\iota_{\widetilde{v}, \widetilde{w}}(\zeta\wedge \omega_X^m)\big|_{F_t}  = \big(\omega_X(\widetilde{v}, \widetilde{w}) \cdot \omega_X^m - \iota_{\widetilde{v}}(\omega_X) \wedge \iota_{\widetilde{w}}(\omega_X^m) + 0 \big)\big|_{F_t}
		\end{align*}
		which is the identity in question.

		\begin{Blank}
			\emph{Step 3: It holds that $\iota_{\widetilde{v}, \widetilde{w}}(\omega_X^{m+1})  =(m+1) \left(\omega_X(\widetilde{v}, \widetilde{w}) \cdot \omega_X^m - \iota_{\widetilde{v}}(\omega_X) \wedge \iota_{\widetilde{w}}(\omega_X^m) \right) $.}
		\end{Blank}
		This is just a straightforward computation:
		\begin{align*}
		\iota_{\widetilde{v}, \widetilde{w}}(\omega_X^{m+1})
		&= (m+1) \cdot \omega_X(\widetilde{v}, \widetilde{w}) \cdot \omega_X^m  - m(m+1) \cdot \iota_{\widetilde{v}}(\omega_X)\wedge \iota_{\widetilde{w}}(\omega_X)\wedge \omega_X^{m-1}\\
		&= (m+1) \cdot \left( \omega_X(\widetilde{v}, \widetilde{w}) \cdot \omega_X^m - \iota_{\widetilde{v}}(\omega_X) \wedge \iota_{\widetilde{w}}(\omega_X^m) \right).
		\end{align*}
		\begin{Blank}
			\emph{Step 4: Conclusion.}
		\end{Blank}
		Combining the results of \emph{Step 2} and \emph{Step 3} we find that
		\begin{align*}
		\iota_{\widetilde{v}, \widetilde{w}}(s(\omega_X)\wedge \omega_X^m)\big|_{F_t} = \frac{1}{m+1} \cdot \iota_{\widetilde{v}, \widetilde{w}}(\omega_X^{m+1})\big|_{F_t}.
		\end{align*}
		Thus, the integrands in \cref{eqs:ClassDefiningCanExt_4} and \cref{eqs:ClassDefiningCanExt_5} above agree (up to scaling) and, hence,
		\begin{align*}
		f_*\left( s(\omega_X)\wedge\omega_X^m\right) (v, w) 
		= \frac{(f_*\omega_X^{m+1})(v, w)}{m+1}, \quad \forall v \in T^{(1,0)}T, \forall w \in T^{(0,1)}T.
		\end{align*}
		This proves \cref{eqs:ClassDefiningCanExt_3} and, as discussed above, the result immediately follows.
	\end{Proof}
	\cref{fibreBundleStructureOnZX} yields a splitting $Z_X \cong Z_{X/T} \times_T Z_T$. We will now concentrate on $Z_{X/T}$:
	
	\begin{Proposition}
		Let $(F, \omega_F)$ be a compact Kähler manifold and consider the complex Lie-group $G := \Aut^0(F)$. Then
		\begin{itemize}
			\item[(1)] the natural action of $G$ on $H^*(F, \mathds{R})$ is trivial.
			\item[(2)] If $H^1(F, \mathds{R}) = 0$, then the action of $G$ on $F$ extends naturally to an action by automorphisms of affine bundles on $Z_{[\omega_F ]}$.
		\end{itemize}
		\label{lifting_action_canonical_extension}
	\end{Proposition}
	
	\begin{Proof}
		Fix an element $g\in G$. As $G$ is connected there exists a path from $\id_F$ to $g$ in $G$. But this is nothing but a homotopy between $\id_F$ and $g$. Thus, all maps in $G$ are null homotopic and so $G$ acts trivially on $H^*(F, \mathds{R})$. This proves $(1)$.
		
		Regarding $(2)$, any element $g \in G$ naturally induces an isomorphism of affine bundles
		\begin{align*}
		g\colon Z_{[\omega_F]} \rightarrow g^*Z_{[\omega_F]} = Z_{[g^*\omega_F]}.
		\end{align*}
		According to item $(1)$, $[g^*\omega_F] = [\omega_F]$ for all $g\in G$ and, hence, there exists \emph{an} isomorphism of affine bundles $Z_{[g^*\omega_F]} \cong Z_{[\omega_F]}$. We claim that in fact there exists only one such isomorphism. In particular, we may functorially identify $Z_{[g^*\omega_F]}$ and $Z_{[\omega_F]}$ and so the action of $G$ on $F$ lifts to $Z_F$ as required.
		
		Regarding the claim, by construction any isomorphism as above is induced by an isomorphism of extensions or, in other words, by a commutative diagram as below: 
		\begin{center}
			\includegraphics{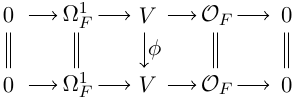}
		\end{center} 
		It is now easily verified by a diagram chase that any morphism $\phi$ making the above diagram commute is of the form $\phi = \id + p \cdot \eta$, where $\eta \in \Hom(\mathcal{O}_F, \Omega_F^1) = H^0(F, \Omega_F^1)$
		and, as before, $V\overset{p}{\rightarrow} \mathcal{O}_X$. But $\dim_\mathds{C} H^0(F, \Omega^1_F) = \dim_\mathds{R} H^1(F, \mathds{R}) = 0$ by the Hodge decomposition. Thus, there is only one isomorphism of affine bundles $Z_{[g^*\omega_F]} \cong Z_{[\omega_F]}$ and we are done.
	\end{Proof}
	
	\begin{Lemma}
		Let $f\colon X \rightarrow T$ be a holomorphic fibre bundle with structure group $G$ and with typical fibre $F$. Suppose that $G\subseteq\Aut^0(F)$ and that $H^1(F, \mathds{C}) = 0$. Then also
		\begin{align*}
		f\circ p \colon Z_{X/T} := Z_{\Omega^1_{X/T}, [\omega_{X/T}]}\rightarrow X  \rightarrow T
		\end{align*}
		is a holomorphic fibre bundle. Its typical fibre is $Z_{F, [\omega_X|_F]}$ and the structure group may be chosen to be $G$.\label{extension=fibre_bundle}
	\end{Lemma}
	Note that $G$ indeed acts on $Z_{F}$ by \cref{lifting_action_canonical_extension} so that the assertion about the structure group of the bundle makes sense.
	
	\begin{Proof}
		Since both $f\colon X\rightarrow T$ and $p\colon Z_{X/T} \rightarrow X$ are holomorphic fibre bundles $f\circ p$ is at least a surjective holomorphic submersion. Moreover, it follows from the functoriality of $Z_{-}$ (see \cite[Lemma 1.16(b)]{greb_CanonicalExtensions}) that the fibre of $f\circ p$ over $t \in T$ is given by
		\begin{align*}
		(f\circ p)^{-1}(t) = p^{-1}(F_t) =  Z_{X/T}\times_X F_t = Z_{\Omega^1_{X/T}|_{F_t}, \left[\omega_X|_{F_t}\right]} = Z_{\Omega^1_{F_t}, [\omega_X|_{F_t}]}.
		\end{align*}
		Now, fix $t\in T$, denote $F := f^{-1}(t)$ and choose a sufficiently small open polydisc $t\in U\subset T$ so that $f^{-1}(U) \cong U\times F$ is trivial. Since $U$ is a polydisc it holds that $H^j(U, \mathds{C}) = 0$ for all $j>0$. Thus, according to the classical Künneth formula the restriction map
		\begin{align*}
		\cdot|_{ \{ t \} \times F}\colon H^*(U\times F, \hspace{0.1cm} \mathds{C}) \rightarrow H^*(F, \mathds{C})
		\end{align*}
		is an isomorphism. In particular, we find that $[\omega_X|_{U\times F}] = pr_{F}^*[\omega_X|_{F}]$. Using again the functionality of extensions and the fact that $\mathcal{T}_{U\times F/U} = pr_F^*\mathcal{T}_F$ we compute
		\begin{align*}
		\left. Z_{\Omega^1_{X/T}, [\omega_{X/T}]}\right|_U 
		= Z_{\Omega^1_{U\times F/U}, [\omega_{X/T}]} 
	    = Z_{pr_F^*\Omega^1_F, pr_F^*[\omega_F]}
        = pr_F^*Z_{F, [\omega_F]} := U \times Z_{F, [\omega_F]}.
		\end{align*}
		This proves that $Z_{X/T} \cong U \times Z_{F}$ as fibre bundles and respecting the affine bundle structure on both sides. We conclude that $f\circ p$ is a holomorphic fibre bundle with fibre $Z_F$.

		The assertion about the structure group being $G$ is clear, because we already saw as part of the proof of \cref{lifting_action_canonical_extension} that given any $g\in G$, there is one and only one identification of $Z_F$ and $g^*Z_F$ as affine bundles. Hence, both $f\colon X\rightarrow T$ and $f\circ p\colon Z_{X/T} \rightarrow T$ are constructed using the same transition functions.
	\end{Proof}
	
	\begin{Remark}
		Record for later reference that both $f\colon X\rightarrow T$ and $f\circ p\colon Z_{X/T} \rightarrow T$ are constructed using the same transition functions. In particular, the first is locally constant if and only if the latter is so.\label{rem:CanExtFibreBunldeFlat}
	\end{Remark}

	\begin{Corollary}
		Let $f\colon X \rightarrow T$ be a holomorphic fibre bundle. Assume that the typical fibre $F$ of $f$ is a Fano manifold. Suppose moreover that the structure group $G$ of $f$ is contained in $\Aut^0(F)$ and that the short exact sequence
		$0 \rightarrow \mathcal{T}_{X/T} \rightarrow \mathcal{T}_{X} \rightarrow f^*\mathcal{T}_{T}
		\rightarrow 0$ splits.
		
		Then there exists an isomorphism of affine bundles
		\begin{align}
		Z_{X, [\omega_X]} \cong Z_{\Omega^1_{X/T}, [\omega_{X/T}]} \times_T Z_{T, [\omega_T]}\label{eqs:SplittingCanExt},
		\end{align}
		where $\omega_T := f_*(\omega_X^{m+1})$. Moreover, the projection map
		\begin{align*}
		\bar{f}\colon Z_{X, [\omega_X]} \rightarrow Z_{T, [\omega_T]}
		\end{align*}
		makes $Z_X$ into a (locally constant if $f$ is so) holomorphic fibre bundle over $Z_T$ with fibre $Z_{F, [\omega_X|_F]}$ and structure group $G$.\label{canExtFibreBundleSTructure}
	\end{Corollary}
	\begin{Proof}
		First of all, \cref{eqs:SplittingCanExt} has already been verified in \cref{fibreBundleStructureOnZX}. Regarding the second assertion, note that $H^1(F, \mathds{C}) = 0$ as $F$ is Fano. Thus, \cref{extension=fibre_bundle} yields that $Z_{X/T} \rightarrow T$ is a (locally constant; see \cref{rem:CanExtFibreBunldeFlat}) holomorphic fibre bundle with structure group $G$ and fibre $Z_F$. But \cref{eqs:SplittingCanExt} just says that 
		\begin{align*}
		\bar{f}\colon Z_{X, [\omega_X]} 
		\rightarrow Z_{T, [\omega_T]}
		\end{align*}
		is the pullback along $Z_T \rightarrow T$ of the bundle $Z_{X/T} \rightarrow T$. Hence, along with $Z_{X/T} \rightarrow T$ also $\bar{f}$ is a (locally constant) holomorphic fibre bundle with structure group $G$ and fibre $Z_F$.
	\end{Proof}
    \begin{Remark}
        In the situation of \cref{canExtFibreBundleSTructure} even if $G$ is not contained in $\Aut^0(F)$ there always exists a finite \'etale cover of $T$ after which we can assume that $G\subseteq \Aut^0(F)$. Indeed, as $F$ is Fano the group $\Aut(F)/\Aut^0(F)$ is finite (cf.\ for example \cite[Corollary 2.17]{brion_AutomorphismGroups}). Moreover, as $G$ acts effectively on $F$ there exists a unique holomorphic principal $G$-bundle $\mathcal{G}\rightarrow T$ such that $X\rightarrow T$ is the associated bundle with typical fibre $F$. Then $T' := \mathcal{G}/G^0\rightarrow T$ is a finite étale cover of $T$ and by construction the structure group of $\mathcal{G}\times_T T'$ may be reduced to $G^0$. In effect, the same is true of the associated bundle $X \times_T T' \rightarrow T'$ and so we are done.\label{covering_trick}
    \end{Remark}
    We are now finally ready to prove \cref{IntroCanExtStein}, the main result of this section:
	
	\begin{Theorem}
		Let $(X, \omega_X)$ be a compact K\"ahler manifold with nef tangent bundle. If the weak Campana-Peternell conjecture \cref{con:Campana-Peternell} holds true then the canonical extension $Z_{X, [\omega_X]}$ is a Stein manifold.
	\end{Theorem}
	\begin{Proof}
		According to \cref{flatness_of_Albanese} there exists a finite \'etale cover $\pi\colon X'\rightarrow X$ such that the Albanese $\alpha\colon X' \rightarrow \Alb(X') =: T$ is a locally constant holomorphic fibre bundle. Its fibres are Fano manifolds with nef and, hence, assuming \cref{con:Campana-Peternell} also big tangent bundle. Possibly replacing $X'$ by another finite étale cover we may assume by \cref{covering_trick} that the structure group $G$ of $\alpha$ is contained in $\Aut^0(F)$. But in this situation \cref{canExtFibreBundleSTructure} applies to $(X', \pi^*\omega_X)$ and shows that there exists a natural map
		\begin{align}
		\bar{\alpha}\colon Z_{\widetilde{X}, [\pi^*\omega_X]} \rightarrow Z_{T, [\omega_T]}\label{eqs:CanExtKahlerNefTangent}
		\end{align}
		making $Z_{X'}$ into a locally constant fibre bundle with structure group $G\subseteq \Aut^0(F)$ and fibre $Z_{F, [\pi^*\omega_X|_F]}$. Here, $\omega_T$ in \cref{eqs:CanExtKahlerNefTangent} above is some (explicitly determined) K\"ahler form on $T$. Note that by \cref{lifting_action_canonical_extension} $\Aut^0(F)$ acts on $Z_F$ so that we may well assume the structure group of $\bar{\alpha}$ to be $\Aut^0(F)$.
		Note moreover, that by the work of \cite[Proposition 2.13]{greb_CanonicalExtensions} $Z_T$ must be Stein as a canonical extension of a complex torus and that due to \cite[Theorem 1.2]{HP_Stein_complements} $Z_F$ is Stein as a canonical extension of a Fano manifold with big and nef tangent bundle. 
		
		In summary, $Z_{X'}$ is naturally a holomorphic fibre bundle over the Stein manifold $Z_T$. The typical fibre of this bundle is $Z_F$, a Stein manifold, and the structure group of the bundle may be chosen to be the connected group $\Aut^0(F)$. But it is a classical theorem by \cite[Théorème 6]{matsushima_SteinFibreBundles} that in this situation also the total space $Z_{X', [\pi^*\omega_X]}$ of the bundle is Stein. Finally, together with $\pi\colon X'\rightarrow X$ also $Z_{X'}\rightarrow Z_{X}$ is a finite étale covering (see \cite[Lemma 2.10.(b)]{greb_CanonicalExtensions}) and we conclude that $Z_{X}$ must be Stein by \cite[Lemma 2]{narasimhan_SteinSpaces}. 
	\end{Proof}

		\section{The special case of surfaces}
	\label{sec:CanExtSurfaces}

    In this section we will provide a proof for \cref{Intro:RuledSurfaces}:

    %Note that item $(3)$ is not quite what we expect: First of all, if $g(C)\geq 2$, then the tangent bundle of $X$ can not be nef and so we would not expect any canonical extension to be Stein. Here, the reason for the first assertion is the relative tangent bundle sequence $0 \rightarrow \mathcal{T}_{X/C} \rightarrow  \mathcal{T}_{X} \rightarrow  \pi^*\mathcal{T}_{C} \rightarrow 0$: If $\mathcal{T}_X$ were nef then so were its quotient $\pi^*\mathcal{T}_C$ and, hence, $\mathcal{T}_C$ itself.
	
	%Moreover, it is well-known that the tangent bundle of a ruled surface $X=\mathds{P}(\mathcal{E})$ over an elliptic curve is nef if and only if the defining bundle $\mathcal{E}$ is semi-stable (cf.\ \cite[Theorem 6.1.]{DPS_ManifoldsWithNefTangentBundle}). This raises the question of what is true in the remaining cases. Indeed, we are able to rule out the higher genus case as well:

	\begin{Proof}[of \cref{Intro:RuledSurfaces}] Assume to the contrary that there exists a Kähler metric $\omega_X$ on $X = \mathds{P}(\mathcal{E})$ whose canonical extension $Z_X$ is Stein. 

    Note that $\pi\colon X \rightarrow C$ is a locally constant fibre bundle as $\mathcal{E}$ is semi-stable; see for example \cite[Theorem 1.5, Proposition 1.7]{jahnke_SemistableRuledSurface} for a proof of this rather basic fact. In other words, if we denote by $\widetilde{C} \overset{p}{\rightarrow} C$ the universal cover of $C$, then there exists a group homomorphism $\rho\colon \pi_1(C) \rightarrow \Aut(\mathds{P}^1)=:G$ such that
		\begin{align*}
		    X \cong \pi_1(C)\backslash(\widetilde{C}\times \mathds{P}^1).
		\end{align*}
		Here, the reason for exceptionally denoting the quotient as one from the left is that shortly we will introduce a second action of a group. It will be crucial below that both of these groups will act from different sides so that the actions commute.
		
		In any case, as $\pi\colon X \rightarrow C$ is a locally constant fibre bundle with fibre $\mathds{P}^1$ - a Fano manifold with connected automorphism group - \cref{canExtFibreBundleSTructure} applies and shows that $Z_X$ is a locally constant fibre bundle over $Z_C$ with typical fibre $Z_{\mathds{P}^1}$ and with the same transition functions as $X\rightarrow C$. Here, for the latter assertion, we use \cref{rem:CanExtFibreBunldeFlat} and the fact, that by \cref{lifting_action_canonical_extension} the action of $\Aut(\mathds{P}^1)$ on $\mathds{P}^1$ lifts uniquely to $Z_{\mathds{P}^1}$. In summary,
		\begin{align}
		Z_{X, [\omega_X]} 
		\cong \pi_1(C)\backslash \left( Z_{\widetilde{C}, [p^*\omega_C]} \times Z_{\mathds{P}^1, [\omega_X|_{\mathds{P}^1}]}\right) 
		= \pi_1(C)\backslash \left( Z_{\widetilde{C}, [p^*\omega_C]} \times G/L \right) .\label{eqs:CanExtSurf1}
		\end{align}
		Here, $\omega_C := f_*(\omega_X\wedge\omega_X)$ is the induced Kähler form on $C$. Moreover, we used that according to \cite[Proposition 2.23]{greb_CanonicalExtensions} there exists a canonical $G$-equivariant identification of canonical extensions
		\begin{align}
		    \left( Z_{\mathds{P}^1} \rightarrow \mathds{P}^1 \right) = \left( G/L \rightarrow G/P \right), \label{rem:CanExtHomFanos}
		\end{align}
        where $L\subsetneq G = \textmd{PGL}_2$ is the maximal diagonal torus and $P\subsetneq \textmd{PGL}_2$ is the Borel subgroup of upper-triangular matrices. In any case, we see that $L\cong \mathds{C}^\times$ is connected and Stein.
		
		Now, let us consider the manifold 
		\begin{align}
		    \mathcal{G}:= \pi_1(C)\backslash \left(  Z_{\widetilde{C}, [p^*\omega_C]} \times G\right).\label{eqs:CanExtSurf2}
		\end{align}
	    The natural projection $\mathcal{G}\rightarrow Z_C$ makes it into a (right) principal $G = \Aut(\mathds{P}^1)$-bundle. Then combining \cref{eqs:CanExtSurf2} with \cref{eqs:CanExtSurf1} we deduce that
		\begin{align*}
		   Z_{X, [\omega_X]} 
		   \cong \pi_1(C) \backslash \left( Z_{\widetilde{C}, [p^*\omega_C]} \times G/L\right)
		   \cong  \pi_1(C) \backslash \left( Z_{\widetilde{C}, [p^*\omega_C]} \times G\right)/L 
		   = \mathcal{G}/L.
		\end{align*}
		In other words, $\mathcal{G} \rightarrow Z_X$ is naturally a (right) principal $L$-bundle. Note that $Z_X$ is Stein by assumption and that $L$ is connected and Stein (cf.\ \cref{rem:CanExtHomFanos}). Therefore, \cite[Théorème 6]{matsushima_SteinFibreBundles} again applies and proves that also $\mathcal{G}$ is Stein. On the other hand, $\mathcal{G}\rightarrow Z_C$ is naturally a (right) $G=\Aut(\mathds{P}^1)$-bundle. Since quotients of Stein spaces by reductive groups are again Stein \cite{snow_reductiveQuotientSteinSpaces}, we infer that also $Z_{C, [\omega_C]} = \mathcal{G}/G$ is Stein. But this contradicts \cite[Example 3.6]{greb_CanonicalExtensions} as $g(C)\geq2$. Thus, $Z_X$ can not be Stein after all and we are done.
	\end{Proof}
    \begin{Remark}
        Note that essentially ad verbatim the same argument also yields the following: Let $f\colon X\rightarrow Y$ be a locallly constant fibration with fibre $F$ and assume that $F= G/P$ is a homogeneous Fano. If the exists a Kähler form $\omega_X$ on $X$ such that the canonical extension $Z_{X, \omega_X}$ is Stein, then there exists a Kähler form $\omega_Y$ on $Y$ so that also $Z_{Y, \omega_Y}$ is Stein.
    \end{Remark}
    Combining \cref{Intro:RuledSurfaces} with \cite[Theorem 1.13]{HP_Stein_complements} this proves \cref{con:GrebWongCanExt} for smooth projective surfaces with the exception of unstably ruled surfaces over elliptic curves.

    % ------------------------------------------------------------------------------------ %											
    %                                   Bibliography									   %
    % ------------------------------------------------------------------------------------ %
	\printbibliography
	
\end{document}